\documentclass[12pt]{amsart}
\usepackage{latexsym,enumerate}
\usepackage{amsmath,amsthm,amsopn,amstext,amscd,amsfonts,amssymb,mathrsfs}
\usepackage{color}
\usepackage{color}
\numberwithin{equation}{section}

\numberwithin{theorem}{section}
\numberwithin{corollary}{section}
\numberwithin{lemma}{section}
\numberwithin{definition}{section}
\numberwithin{proposition}{section}
\numberwithin{remark}{section}

\begin{document}

\title[Symmetry for a general class of overdetermined elliptic problems]
{Symmetry for a general class of overdetermined  elliptic problems}

\author{ F. Brock$^1$ }
\thanks{}
\date{}

\begin{abstract}
Let $\Omega $ be a bounded domain in $\mathbb{R} ^N $, and 
let $u\in C^1  (\overline{\Omega }) 
$ be a weak solution of the following 
overdetermined BVP: $-\nabla (g(|\nabla u|)|\nabla u|^{-1} \nabla u 
)=f(|x|,u)$, $ u>0 $ in $\Omega $ and  $u(x)=0, \ |\nabla u (x)| =\lambda 
(|x|)$ on 
$\partial \Omega $, where $g\in C([0,+\infty  ))\cap C^1 ((0,+\infty ) ) $ with 
$g(0)=0$, $g'(t)>0$ for $t>0$, $f\in C([0,+\infty )) \times [0, +\infty ) )$, 
$f$ is nonincreasing in $|x|$, $\lambda \in C([0, +\infty )) $ and $\lambda $ 
is positive and nondecreasing. We show that $\Omega $ is 
a ball and $u$ satisfies some "local" kind of symmetry. The proof is 
based on the method of continuous Steiner symmetrization.  
\noindent 
 \bigskip

\textsl{Key words:} degenerate elliptic equation, overdetermined 
boundary value problem,
symmetry of the solution, continuous rearrangement 

\textsl{2000 Mathematics Subject Classification:} 28D10, 35B05, 35B50, 35J25, 35J60, 35J65
\end{abstract}

\maketitle

\setcounter{footnote}{1} 
\footnotetext{F. Brock: University of Rostock, Department of 
Mathematics, 
Ulmenstr. 69, 18057 Rostock, Germany, email: friedemann.brock@uni-rostock.de}

\section{ Introduction}
In a celebrated paper \cite{Se} Serrin proved the following symmetry 
result for an overdetermined elliptic boundary value problem:
\\[0.1cm]
{\bf Theorem A.} {\sl Let $\Omega \subset \mathbb{R}^N $ be a 
bounded domain with $C^2 $-boundary, and let $u\in C^2 (\overline{\Omega })$ 
satisfy
\begin{eqnarray}
 & & Lu\equiv a(u,|\nabla u|) \Delta u
 + \sum_{i,j=1}^N b(u,|\nabla u|) u_{x_{i}} u_{x_{j}} u_{x_{i} x_{j} }
=c(u,|\nabla u|) ,
\nonumber 
\\
\label{serrin1}
 & &  u>0\quad \mbox{in } \ \Omega , 
 \\
\label{serrin2}
 & & u=0,\ |\nabla u| = constant  \quad \mbox{on } \ \partial \Omega ,
\end{eqnarray}
where $a,b$ and $c$ are continuously differentiable in each variable 
and $L$ is uniformly elliptic. Then 
$\Omega $ is a ball and $u$ is radially symmetric about the center of 
the ball.
\/}
\\[0.1cm]
The proof of Theorem A in \cite{Se} uses the so-called {\sl moving plane 
method\/} which became very popular thanks to Gidas, Ni and 
Nirenberg's paper \cite{GiNiNi}. The method combines symmetry arguments and 
boundary versions of the strong maximum principle, and it 
has often been applied to show the symmetry of solutions in overdetermined 
problems, e.g. in 
\cite{AfBu3}, \cite{Co}, \cite{Dal}, \cite{KuPr}, 
\cite{PhRa}, 
\cite{Re1}-\cite{Re3}, 
\cite{Vo}, \cite{WiGlSi}. The moving plane device applies 
to very general - even fully nonlinear - elliptic equations. On the 
other hand, if the equation degenerates and/or contains terms which are 
less regular, then the method often fails. 
\\ 
Several other tools have been applied in such situations. Let us give a short overview.    
One approach is based on a comparison principle   
which is combined with some Rellich-type identity, 
see e.g. \cite{We}, \cite{GaLe}, 
\cite{Mo}, 
\cite{Ra2}, \cite{FraGaKa} and \cite{FaKa}.   
Another idea, used e.g. in \cite{Be}, \cite{PaSc2}, is to exploit some integral 
identity which is equivalent to the overdetermined problem. 
Although these two methods can be applied to 
degenerate operators - for instance, to the $p$-Laplacian - they are useful only 
for very special equations. 
A third method is based on a comparison with suitable radial solutions of the equation, 
and it is applicable to situations when the solution of the bondary value problem is unique, see e.g.
\cite{HePhPr}, 
\cite{Gre1} and \cite{Gre2}.
A fourth approach is based on 
the {\sl method of domain derivative\/} which has been widely investigated 
in shape optimization (see \cite{SoZo}). This device again seems useful in 
problems where the solution of the boundary value problem is unique
(see \cite{ChoHe}, \cite{HePh}). 
Note, this approach also highlights the relation between a second ('overdetermined') boundary condition and minimization of appropriate domain functionals, see \cite{BaWa1},\cite{BaWa2}. The method of domain derivative has also been combined with 
another tool: the so-called {\sl continuous Steiner 
symmetrization\/} (CStS) (see \cite{ChaChoHe}, \cite{BrHe}). The idea of CStS is 
to find "local analogues" to some well-known rearrangement 
inequalities  (see \cite{PoSz}, \cite{Kaw2}, \cite{Br1}, \cite{Br2}). The author 
exploited this method to prove symmetry results 
for nonnegative solutions of boundary value 
problems in symmetric domains (see \cite{Br2}, \cite{Br3}).  \\  
The aim of this paper is to give a new approach to overdetermined 
problems which is  based on the CStS, but does not use domain derivatives.  
Although our method is restricted to operators in  divergence form, 
we allow nonsmooth terms in the equation, and the 
solution of the boundary value problem need not be unique. 
\\
We fix some 
notation. By $x=(x_{1} ,\ldots ,x_{N} )$ we denote a point in 
$\mathbb{R}^N $, and by $|x|$ its norm. Our main result is
\\[0.1cm]
{\bf Theorem 1:} {\sl
Let
$f: [0, +\infty )\times [0, +\infty ) \to \mathbb{R} $ 
be a bounded measurable function, such that
the mapping 
$ v\mapsto f(r,v)$
 is continuous, uniformly for all $r$, and the mapping 
$r\mapsto f(r,v)$ is nonincreasing, ($r,v\in [0, +\infty )$). 
Let $g\in C([0,+\infty) )\cap C^1 ((0,+\infty ) ) $, with 
$g(0)=0$, $g'(t)>0$ for $t>0$, and let
$\lambda \in C([0, +\infty )) $ be a positive
and nondecreasing function.
Further, let $\Omega $ be a bounded domain, and let   
$u\in C^1 (\Omega ) \cap C (\overline{\Omega})$ be a 
weak solution of the following problem,
\begin{eqnarray}
\label{eq1}
 & & -\nabla \left( g(|\nabla u|)\frac{\nabla u}{|\nabla u|} \right) 
=f(|x|, u),\ u>0 \quad \mbox{in } \ \Omega ,
\\
\label{bdry1}
 & & u(x) =0 \quad \mbox{on } \ \partial \Omega , 
\\  
\nonumber
 & & \mbox{ Given $\varepsilon >0$, there is an open set $U$ containing $\partial \Omega $ such that}
\\
 & &  
\label{bdry2nhood}
\left| |\nabla u(x) |-\lambda (|x |) \right| <\varepsilon \quad \forall x\in U\cap \Omega .
\end{eqnarray}
Then $\Omega $ and $u$ 
satisfy the following symmetry properties: 
\\
$\Omega $ is ball,
\begin{eqnarray}
\label{decomp1}  
 & & \Omega =  \left( \bigcup\limits_{k=1}^{m} C_{k }\right) \cup S, 
 \quad \mbox{ where} 
\\ 
 & & C_{k } := \left\{ x\in \Omega :\, r_k <|x-y^k | <R_k \right\} 
 ,\quad  R_k >r_k \geq 0, \ y^k \in \Omega  ,
\nonumber 
\\  
 & & \frac{\partial u}{\partial \rho }   <  0 \quad  \mbox{ in } \ 
 C_{k} ,\qquad ( \rho  := |x-y^{k} |  ),  
\nonumber 
\\  
 & & u(x) \geq u(y) \quad  \mbox{if } \  0\leq |x-y^k | \leq |y-y^k | =r_k , \qquad 
1\leq k\leq m, 
\qquad \mbox{and}
	\nonumber 
\\    
 & & \nabla u = 0 \qquad \mbox{in } \ S. 
\nonumber
\end{eqnarray}  
The sets on the right-hand side of (\ref{decomp1}) are mutually disjoint and        
there can be a countable number of annuli $C_{k} $, i.e.   
$m=+\infty $. }
\\[0.1cm]
{\bf Remark 1:}  
\\
{\bf (a)} 
(\ref{eq1}) means 
\begin{equation}
\label{eq1weak}
\int\limits_{\Omega } 
g(|\nabla u|) \frac{\nabla u \cdot  \nabla v }{| \nabla u|} \, dx = \int\limits_{\Omega }
f(|x|, u)v \, dx \qquad \forall v\in W_0 ^{1,1} (\Omega )\cap W^{1, \infty } (\Omega ), 
\end{equation}
where
the expression $g(|y|) \frac{y}{|y|} $, ($y\in \mathbb{R} ^N $),  
is interpreted as the zero vector, if $y=0$. 
\\[0.1cm] 
Note that if
$\partial \Omega $ is smooth and $u\in C^1 (\overline{\Omega })$, then 
(\ref{bdry2nhood}) means 
\begin{equation}
\label{bdry2} 
\frac{\partial u}
{\partial \nu} (x) = \lambda (|x|) \ \mbox{ on } \partial \Omega , \qquad (\nu: \ 
\mbox{ exterior unit normal }).
\end{equation}
{\bf (b)} Theorem 1 falls out of the scope of the above mentioned results for the following reasons: 
\\ 
${\bf \cdot }$  We do not assume that $f$ is smooth in the second variable.
\\
${\bf \cdot }$ The solution of the boundary value problem (\ref{eq1}), (\ref{bdry1}) might be not unique.   
\\
${\bf \cdot }$ The differential operator in (\ref{eq1}) is not assumed to be uniformly elliptic.
\\  
We also emphasize that the solution $u$ is not radially symmetric 
in general. For instance, there are examples of nonsymmetric 
solutions of problem (\ref{eq1}),(\ref{bdry1}) in a ball in the 
$p$-Laplacian case, that is if  
$g(z) =z^{p-1} $ for some $p>1$ (see \cite{Br2}, \cite{BrHe}).
\\[0.1cm]  
Now we outline the content of the article. 
In section 2, we give the definition of CStS and we present some 
results of \cite{Br1}-\cite{Br3} which will be of later use.   
In Section 3 we prove Theorem 1 and we give some extensions of 
the result in Theorem 2.
\section{ Continuous Steiner symmetrization}  
For points $x\in \mathbb{R}^N $, ($N\geq 2$), we write $x=(x_{1} 
,x')$, where 
$x'=(x_2 ,\ldots ,x_{N} )$, and for any set $M\subset \mathbb{R} ^N $ let $\chi (M)$ 
the characteristic function of $M$. If $u:\mathbb{R}^N \to \mathbb{R}$, 
then let $\{ u>a \} $ and $\{ b\geq u>a\} $ denote the sets
$\{ x\in \mathbb{R} ^N :\, u(x)>a\} $, and $\{ x\in \mathbb{R}^N :\, 
b\geq u(x)>a\} $, respectively, ($a, b\in \mathbb{R} $, $a<b $).       
Let $\mathscr{L} ^k$ denote $k$-dimensional Lebesgue 
measure, ($1\leq k\leq N$), and $\Vert \cdot \Vert _{p} $ 
the usual norm in $L^p (\mathbb{R} ^N )$, ($1\leq p\leq +\infty $). By 
$\mathscr{M} (\mathbb{R}^N )$ we denote the family of   
Lebesgue measurable  - measurable in short -  
sets in $\mathbb{R}^N $ with {\sl finite\/} measure. Finally, 
let $\mathscr{S} _+ (\mathbb{R} ^N )$  denote the class of real, 
nonnegative  measurable   
functions $u$ satisfying  
$$   
\mathscr{L}^N \left( \{   
 u>c\} \right) <+\infty \quad \forall c>0.  
$$  
Note that nonnegative 
functions in $L^p (\mathbb{R}^N )$, ($1\leq p< +\infty $), belong to
 $\mathscr{S} _+ (\mathbb{R}^N )$. Generally we  treat measurable 
sets and functions in {\sl a.e. sense\/}.
\\   
Given a unit vector $e\in \mathbb{R} 
^n $, a {\sl continuous Steiner symmetrization\/} (CStS) 
is a continuous homotopy 
which connects sets $M\in \mathscr{M} (\mathbb{R} ^N ) $ and functions $u\in 
\mathscr{S} _{+} (\mathbb{R} ^N )$ with their  Steiner symmetrizations in 
direction $e$, $M^* $, 
respectively $u^{*} $. Homotopies of such type can be constructed in different  
ways (see \cite{Br1}, \cite{Br2} and the references cited therein). Below we 
define a variant of CStS which has been investigated by the author in 
\cite{Br1}, \cite{Br2}. 
\\[0.1cm]
For the convenience of the reader we 
first recall the definition of the   
well-known Steiner symmetrization (see e.g. \cite{Kaw1}).   
\\[0.1cm]
{\bf Definition 1:} {\sl (Steiner symmetrization)
\\  
{\bf (i)} For any set  
$M\in \mathscr{M} (\mathbb{R} )$ let   
$$  
M^* := \left( -\frac{1}{2} \mathscr{L}^1 (M),\frac{1}{2} 
\mathscr{L}^1 (M) \right)   
$$   
denote the (one-dimensional) symmetrization of $M$.
\\   
{\bf (ii)} Let   
$M\in \mathscr{M } (\mathbb{R}^N )$, ($N\geq 2$). For every $x'\in \mathbb{R}^{N-1} $   
let  
$$  
M(x') := \left\{ x_{1}\in \mathbb{R} :\ (x_{1} ,x')\in M \right\} .  
$$  
The set   
\begin{equation}  
M^*  := \left\{ x=(x_{1},x'):\ x_{1} \in \left( M(x')\right) ^* ,\   
x'\in \mathbb{R}^{N-1} \right\}   
\label{St1}  
\end{equation}   
is called the {\sl Steiner symmetrization\/} of $M$ (with respect to 
$x_{1}$).
\\ 
(Note that $M^* $ is symmetric and convex with respect to the hyperplane   
$\{ x_{1}=0\} $.)
\\  
{\bf (iii)} If $u\in \mathscr{S} _+ (\mathbb{R}^N )$, ($N\geq 2$),  then the function 
\begin{equation} 
\label{St2}
u^*(x):= 
\left\{ 
\begin{array}{ll}
\sup \left\{ c>0:\ x \in \{ u>c\} ^* \right\} 
 & \ \mbox{ if } \ x \in 
\bigcup\limits_{c>0 } \{ u >c\} ^{*}  \\
0 & \ \mbox{ if } \ x\notin 
\bigcup\limits_{c>0 } \{ u>c\} ^*  
\end{array} 
\right. ,
\quad x\in 
\mathbb{R} ^N ,
\end{equation}  
is called the {\sl Steiner symmetrization\/} of $u$ (with respect to 
$x_{1} $). 
\\  
(Note that $u^* (x_{1} ,x') $ 
is symmetric with respect to $\{ x_{1} =0\} $ and nonincreasing  
in $x_{1}$ for $x_{1}>0$.)}
\\[0.1cm]  
{\bf Definition 2:} {\sl   
(Continuous symmetrization of sets  
 in $\mathscr{M}(\mathbb{R})$)
\\  
 A family of set transformations  
$$  
E_{t } :\quad \mathscr{M }(\mathbb{R} ) \longrightarrow \mathscr{M } (\mathbb{R} ),   
\qquad 0\leq t \leq +\infty ,  
$$  
satisfying the properties,  
$(M,N\in \mathscr{M }(\mathbb{R}),\ 0\leq s,t\leq +\infty )$  
\\
{\bf (i)} 
$\mathscr{L}^1 ( E_{t } (M))  = \mathscr{L} ^1 (M)  ,\quad $   
  (equimeasurability),  
\\
{\bf (ii)} if $M\subset N$, then $E_{t } (M)\subset E_{t } (N),  
\quad $ (monotonicity),
\\  
{\bf (iii)} $E_{t } (E_{s } (M)) =E_{s+t } (M),\quad $   
(semigroup property),
\\  
{\bf (iv)} if $M$ is an interval $[x-R,x+R]$, ($x\in \mathbb{R} $, 
$R>0$),  then \\
$E_{t} (M) :=  
[xe^{-t} -R, xe^{-t} +R ]$,
\\  
is called  continuous symmetrization.}
\\[0.1cm]   
The existence and  uniqueness of the family   
$ E_{t } $, $  
0\leq t \leq +\infty $, has been proved in \cite{Br2}, Theorem 2.1.
\\[0.1cm]
{\bf Definition 3:}  
{\sl ( Continuous Steiner symmetrization (CStS) )  
\\  
{\bf (i)} Let $M\in \mathscr{M} (\mathbb{R}^N)$, ($N\geq 2$). The family of sets  
\begin{equation}  
E_t (M) := \left\{ x=(x_{1} ,x'):\ x_{1} 
\in E_t \left( M(x')\right) ,\ x'\in \mathbb{R}^{N-1}
\right\} ,\    0\leq t \leq +\infty ,  
\label{CSt1}  
\end{equation}  
is called the continuous Steiner symmetrization (CStS) of $M$ (with respect 
to $x_{1}$).
\\  
{\bf (ii)} Let $u\in \mathscr{S}_{+} (\mathbb{R} ^N )$. 
The family of functions $E_t (u) $, $0\leq t \leq +\infty ,$   
defined by  
\begin{equation} 
\label{CSt2}
E_t (u)( x):= 
\left\{ 
\begin{array}{ll}
\sup \left\{ c>0:\ x\in E_t (\{ u>c\} ) \right\} 
 & \ \mbox{ if } \ x\in 
\bigcup\limits_{c>0 } E_t (\{ u>c\} )  \\
0 & \ \mbox{ if } \ x\notin 
\bigcup\limits_{c>0 } E_t (\{ u>c\} )  
\end{array} 
\right. ,
\quad x\in 
\mathbb{R} ^N ,
\end{equation} 
is called  CStS of $u$ with 
respect to $x_{1}$ in the 
case $N\geq 2 $ and continuous
symmetrization in the case $N=1$.}
\\[0.1cm]
{\bf Remark 2. }
\\
{\bf 1.} For convenience, we will henceforth simply 
write $M^t $ and $u^t $ for the sets $E_t (M) $, respectively for the functions $E_t (u) $, ($t\in [0,+\infty ]$).
\\
{\bf 2.}       
It can be shown that, if $M\in \mathscr{M} (\mathbb{R} ^N )$
and $M$ is open, then the sets $M^t  $,
($t\in [0,+\infty ]$),  have open representatives.  
This makes it possible to give  {\sl pointwise\/} definitions of {\sl
open sets\/} and of {\sl continuous functions\/}:
\\[0.1cm]   
{\bf (i)} If $M\in { \mathcal M } (\mathbb{R} )$ is open and $t\in [0,+\infty ]$, then let
\begin{equation}
\label{CSt3}
M^{t,O}  := \bigcup \left\{ U: \ U\mbox{ is an open representative of } N^t  ,\ N 
\mbox{ open },\
 N\subset \subset M \right\} .
\end{equation}
{\bf (ii)} If $M\in { \mathcal M } (\mathbb{R} ^N)$, ($N\geq 2$), is open and $t\in 
 [0,+\infty ]$, then let 
\begin{equation}
\label{CSt4}
M^{t,O }:= \left\{ x=(x_{1} ,x'):\ x_{1}\in  (M(x'))^{t,O } ,\ x'\in \mathbb{R} 
^{N-1} \right\} . 
\end{equation}
Note that the relations (\ref {CSt3}), (\ref{CSt4}) have to be 
understood in {\sl pointwise\/} sense. The sets $M^{t,O} $ in 
(\ref{CSt3}) and (\ref{CSt4}) are open and they are called the 
{\sl precise representatives of\/} $M^t $. 
\\[0.1cm]
{\bf (iii)}
If $u\in \mathscr{S} _{+} (\mathbb{R}  ^N )$ is  continuous and 
$t\in [0,+\infty ]$, then there exists a unique continuous 
representative of $u^t $ which is given by (\ref{CSt2}) - now in 
pointwise sense! -  
where the sets $\{ u^t >c\}  $ have to be replaced by their precise 
representatives.
\\[0.1cm]
From now on let us agree that, if we speak about the CStS of open sets 
or continuous functions then we always mean their {\sl precise\/} 
representatives.
\\[0.1cm]
{\bf Remark 3.}  Below we
summarize basic properties of CStS, which have been proved by the author 
in \cite{Br1}, \cite{Br2},    
($ M\in \mathscr{M} (\mathbb{R}^N ), 
u,v \in \mathscr{S} _+ (\mathbb{R} ^N ) ,\ t \in [0,+\infty ] $).  
\\  
{\bf 1.} {\sl Equimeasurability:\/} From Definitions 2 and 3 we have   
\begin{equation}   
\mathscr{L} ^N (M)=\mathscr{L} ^N (M^{t } ) 
\quad \mbox{and} \quad \{ u^{t } >c\}  =  \{ u>c\} ^{t }
\quad
\forall c>0.  
\label{equimeas}  
\end{equation}
{\bf 2.} {\sl Monotonicity:\/} If $u\leq v$ then $u^t \leq v^t $.
\\
{\bf 3.} If $\psi : [0, +\infty )  \rightarrow [0, +\infty )  $ is bounded 
and nondecreasing with $\psi (0)=0$, then 
\begin{equation}
\label{monot}
\psi (u^t ) = \left( \psi (u)\right) ^t .
\end{equation} 
{\bf 4.} {\sl Homotopy:\/}  We have 
\begin{equation}
\label{homot}
M^0 =M, \   
u^0 =u ,\ M^{\infty } =M^{*} , \ u^{\infty } =u^* .
\end{equation}
Furthermore, if $M=M^{*} $ or $u=u^{*} $, then we have $M=M^t $, 
respectively $u=u^t $. Finally, 
if $t _n \to t $ as $n\to +\infty $ and 
$u\in L^p (\mathbb{R} ^N )$ for some $p\in [1,+\infty )$, then
\begin{equation}  
\lim_{n\to \infty } \Vert u^{t _n } -u^t \Vert _{p} =0.  
\label{cont}  
\end{equation}    
{\bf 5.} {\sl Cavalieri's principle:\/}  
If $F$ is continuous and if $F(u) \in L^1 (\mathbb{R} ^N )$ then 
\begin{equation}  
\int\limits_{\mathbb{R} ^N } F(u)\, dx =\int\limits_{\mathbb{R} ^N } F(u^{t } )\, dx.  
\label{Cav}  
\end{equation} 
{\bf 6.} {\sl Nonexpansivity in $L^p $:\/} 
If $u,v\in L^p (\mathbb{R} ^N ) $ for some $p\in [1,+\infty )$ then
\begin{equation}
\label{nonexp}
\Vert u^t -v^t \Vert _{p} \leq \Vert u-v\Vert _{p } .
\end{equation}
{\bf 7.} {\sl Hardy-Littlewood inequality:\/} 
If $u,v \in L^2 (\mathbb{R} ^N ) $ then 
\begin{equation}
\label{hardylit}
\int\limits_{\mathbb{R} ^N } u^t v^t \, dx \geq \int\limits_{\mathbb{R} ^N } uv \, dx.
\end{equation}
{\bf 8.} A generalization of {\bf 5.} and {\bf 7.} is the 
following result:
\\  
Let $u\in L^{\infty } (\mathbb{R} ^N )$ and suppose that $u$ vanishes outside   
some ball $B_R ,\ R>0 $. Furthermore, suppose that $F=F(x,v)$ is bounded 
and measurable on $B_R \times [0, +\infty ) $, 
continuously differentiable in $v$ with $F(x,0)=0 $ $\ \forall x\in 
\mathbb{R} ^N $,  and $ (\partial /\partial v) F(x,v) $  is even in $x_{1}$ and 
nonincreasing in $x_{1}$ for $x_{1} >0$.  
Then 
\begin{equation}  
\int_{B_R } F(x,u)\ dx \leq \int_{B_R } F(x,u^{t })\ dx .     
\label{FuFut}  
\end{equation}  
{\bf 9.} If $u$ is Lipschitz continuous   
with Lipschitz constant $L$, then $u^{t } $ is Lipschitz continuous, 
too,   
with Lipschitz constant less than or equal to $L$. 
\\
{\bf 10.} If $\mbox{supp } u \subset B_{R}$ for some $R>0$, then we also have  
$\mbox{supp } u^t \subset B_{R} $. 
If, in addition, $u$ is Lipschitz continuous on $\mathbb{R}^N $ with Lipschitz constant 
$L$, then we have 
\begin{eqnarray}
\label{ut-u}
 & & | u^t (x) -u(x) | \leq LRt  \quad \forall x\in B_{R} \quad  
 \mbox{and } \\
 & &   
\int\limits_{B_{R} } G(|\nabla u^{t } | )\ dx \leq 
\int\limits_{B_{R} } G(|\nabla u| )\ dx ,   
\label{polyasz1}  
\end{eqnarray}
for every convex function $G:[0,+\infty )  \rightarrow [0, +\infty ) $ with 
$G(0)=0$. 
\\[0.1cm]  
Note that {\bf 1.}-{\bf 3.}, {\bf 5.}-{\bf 7.}, {\bf 9.} 
and (\ref{polyasz1})  are common 
properties of many rearrangements (see \cite{Kaw1}).
\\   
The following symmetry criteria have been proved in \cite{Br2}, section 6.
\\[0.1cm]
{\bf Lemma 1:} {\sl (see \cite{Br2}), Theorem 6.2)   
Let $\Omega $ be a bounded open set, 
$u\in C^1 (\Omega ) \cap C (\overline{\Omega }) $,  
 $u>0$ in $\Omega $ and  $u=0$ on $\partial \Omega $. Furthermore,  
let $G:[0, +\infty )   \rightarrow [0, +\infty )$ be strictly convex 
with $G(0)=0$, and suppose that  
\begin{equation}  
\lim_{t \to 0 } \frac{1}{t } \Bigg( 
\int\limits_{ \Omega  } G(|\nabla u|)\, dx -  
\int\limits_{ \Omega ^{t }  } G(|\nabla u^{t } |)\ dx \Bigg)  =0.  
\label{limG}  
\end{equation}  
Then $u$ satisfies the following symmetry property:
\\ 
If  $y =(y_1 ,y' )\in \mathbb{R} ^N $ with  
\begin{equation}  
0<u(y )< \sup u,\quad \frac{\partial u}{\partial x_{1}} (y)>0,   
\label{sym1}  
\end{equation}  
and $\widetilde{y} $ is the (unique) point satisfying  
\begin{equation}  
\widetilde{y} =
(\widetilde{y_{1}},y' ),\quad  \widetilde{y_{1}} >y_1 ,
\qquad u(y )=u(\widetilde{y} )<u(z,y') 
\quad    \forall z\in (y_1 ,\widetilde{y_1} ),  
\label{sym2}  
\end{equation}  
then  
\begin{eqnarray}  
\frac{\partial u}{\partial x_i } (y ) 
 & = & \frac{\partial u}{\partial x_i } 
(\widetilde{y} )  ,\quad i=1,\ldots ,N-1, \quad \mbox{ and} 
\nonumber 
\\  
\frac{\partial u}{\partial x_{1}} (y ) & = & -  
\frac{\partial u}{\partial x_{1}} (\widetilde{y} ).   
\label{sym3}  
\end{eqnarray}
} 
We will say that $u$ {\sl is locally symmetric in direction\/} 
$x_{1} $ if $u$ satisfies the properties 
(\ref{sym1})-(\ref{sym3}).
\\
{\bf Lemma 2:} {\sl (see \cite{Br2}, Theorem 6.1) Let $\Omega $, $u$ and $G$ be as in Lemma 1, and suppose that 
for arbitrary  rotations $x\longmapsto y =(y_{1} ,y') $ of  the 
coordinate system, $u$ is locally symmetric in direction $y_{1} $. 
Then $\Omega $ is an at most countable union of mutually disjoint open balls and 
$u$ satisfies the symmetry properties 
(\ref{decomp1}).
}
\section{ Symmetry of the solution} 
In this section we show Theorem 1. The idea is to  
use appropriate test functions $v$ in (\ref{eq1weak}) 
and then to use Lemma 1. This works well in the case of 
{\sl Steiner symmetric\/} domains $\Omega $ (see \cite{Br2}), choosing $v= u^t $, respectively  $v=u$, 
($u^t$: CStS of $u$, with small $t>0$). However, in our situation, 
$\Omega $ is not assumed to be symmetric, so that $u^t $ might not vanish on $\partial \Omega $. 
Therefore we modify the approach of \cite{Br2}, using appropriate cut-off functions of $u^t $ and $u$.      
\\ 
First we introduce some notation. For functions $v$ we 
write $v_{+} := \max \{ v, 0\} $. By the symbol 
$o(t) $ we denote any function satisfying 
$\lim _{t\to 0}o(t)/t =0 $ and which may vary from line to 
line. For any point $x\in \overline{\Omega }$ we write
$$
d(x):= \mbox{dist}\, \{ x; \partial \Omega \} 
\equiv \inf \{ |x-z| :\, z\in \partial \Omega \} .
$$ 
Throughout this section, let $u$ 
be the solution of problem (\ref{eq1})-(\ref{bdry2}).  
For convenience, we extend $u$ by zero 
outside $\Omega $, so that $u\in C^{0,1} (\mathbb{R} ^N )$. 
We denote by $L$ the Lipschitz constant of $u$, and we set $u_0 := \max \{ u(x):\, x\in \Omega \} $.  
We choose a number $R>0$ such that $\overline{\Omega } \subset B_R 
$, and we set
$f_0 := \sup \left\{ |f(|x|,v)|:\, |x|\leq R,\, 0\leq v\leq u_0 \right\} $ and   
$$
k:= 2RL.  
$$   
Finally,  we fix some coordinate 
system 
$$
x=(x_{1} ,x'), \quad  ( x_1 \in \mathbb{R} , \, x' \in \mathbb{R} ^{N-1} ). 
$$
Let $u^t $, ($ 0\leq t\leq +\infty $),  denote the 
CStS of $u$ with respect to $x_{1} $. Since $u\in C^{0,1} (\mathbb{R} ^N )$,  
we have by Remark 2, {\bf 9.} and by (\ref{ut-u}),
\begin{eqnarray}
\label{utLip} 
 & & u^t \in C^{0,1} (\mathbb{R}^N ),   
\\
\label{utL}
 & & u^t \ \mbox{ has Lipschitz constant less than or equal to $L$, and}
\\ 
\label{ut-u0}
 & & |u^t (x)-u(x)|  \leq  LRt 
\qquad \forall t\in [0,+\infty  ],
\end{eqnarray}  
Next we obtain some estimates for $u$ and $u^t $ near the boundary of $\Omega $.
\\
Since $u $ is positive in $\Omega $ and continuous on $\mathbb{R} ^N $, 
and since $\Omega  $ is bounded, we have
\begin{equation} 
\label{udx}
\lim\limits_{s\to 0} \sup \left\{ d(x) : 0<u(x) \leq s \right\} =0.
\end{equation} 
{\bf Lemma 3:} {\sl Let $t\in [0, +\infty )$ and $u^t (x) > kt$. 
Then $x\in \Omega $.}
\\[0.1cm]
{\sl Proof:}  We have by property (\ref{ut-u0}),
\begin{eqnarray*}
u(x) & \geq  & = - |u(x) -u^t (x)| + u^t (x) 
\\
 & \geq  & -RL t + kt =RLt >0.
\end{eqnarray*}
Hence $x\in \Omega $.
$\hfill \Box $
\\[0.1cm]
For convenience, we set
\begin{eqnarray}
\label{At} 
M_1 (t) & := & \{ 0<u \leq kt\} \quad \mbox{and }\\
M_2 (t) & := & \Omega \cap \{ 0<u^t \leq kt \} , \qquad (t\in (0, +\infty) ).
\end{eqnarray}
Note that $M_1(t) , M_2 (t) \subset \Omega $, and by Remark 3, {\bf 1.},
\begin{equation}
\label{M1=M2}
\mathscr{L} ^N (M_1 (t)) = \mathscr{L} ^N (M_2 (t)) \quad \forall t\in (0, +\infty ).
\end{equation}
{\bf Lemma 4:} {\sl There holds: 
\begin{equation}
\label{dist1}
\lim\limits_{t\to 0}  \sup \left\{ d(x) :\, x\in M_1 (t) \cup M_2 (t)  \right\} =0.
\end{equation}  
}
{\sl Proof: } In view of property (\ref{udx}) it is sufficient to show that
\begin{equation}
\label{dist2}
\lim\limits_{t\to 0}  \sup \left\{ d(x) :\, x\in  M_2 (t)  \right\} =0.
\end{equation}  
Assume that (\ref{dist2}) is not true. Then there exists a number $\delta >0$, 
a sequence of points $\{ x_n \} \subset \Omega $ and  a decreasing sequence 
$\{ t_n \} \subset \mathbb{R}$ 
with $\lim\limits_{n\to \infty } t_n =0$, such that 
$0<u^{t_n } (x_n ) \leq kt_n $, but $d(x_n ) \geq \delta $. The latter also implies that $u(x_n ) \geq \varepsilon $, for some $\varepsilon >0$.  On the other hand, we have by (\ref{ut-u}),  
\begin{eqnarray*}
u(x_n ) & \leq & | u(x_n ) - u^{t_n } (x_n ) |  + u^{t_n } (x_n ) 
\\
 & \leq & LR t_n + kt_n \to 0 ,
\end{eqnarray*}
as $n\to \infty $,
a contradiction.  
$\hfill \Box $ 
\\[0.1cm]
{\bf Lemma 5:} {\sl There exists a constant $c_0 >0$ such that 
\begin{equation}
\label{estu<kt}
\mathscr{L}^N ( M_1 (t)  ) \leq c_0 t, \quad (0<t< +\infty ).
\end{equation}
}
{\sl Proof: } By (\ref{bdry2nhood}) and (\ref{udx}) 
there exist positive numbers $\tau$ and $t_0 $, such that     
\begin{equation}
\label{nablaubelow}
|\nabla u(x)| \geq \tau \quad \mbox{if $x\in M_1 (t_0 )$.}
\end{equation}
By the Implicit Function Theorem, we have that for every $s\in (0, kt_0 ]$, 
$\{ u>s \} $ is an open subset of $\Omega $ with $\partial \{ u>s \} = \{ u=s\} $, and $\{ u=s \} $ is locally a $C^1 $--hypersurface.
Integrating (\ref{eq1}) over $\{ u>s \} $, ($s\in (0, k t_0 ]$), Green's Theorem yields
$$
\int\limits_{\{ u=s\} } g(|\nabla u|)\, d\mathscr{H}_{N-1} (x) =  
\int\limits_{\{ u>s\} } f(|x|, u)\, dx. 
$$
By (\ref{nablaubelow}) this implies 
\begin{equation}
\label{bddper}
\int\limits_{\{ u=s\} } \, d\mathscr{H}_{N-1} (x) \leq  
\frac{ f_0 }{g(\tau)}  \mathscr{L} ^N (\Omega ), \quad (s\in (0, kt_0 ] ).   
\end{equation}
Using this and the co-area formula (see \cite{Fe}), we obtain:
\begin{eqnarray*}
\mathscr{L} (M_1 (t) )
 & = &
\int\limits_0 ^{kt} \left( \int\limits_{\{ u=s \} } 
\frac{d\mathscr{H} _{N-1} (x)}{|\nabla u|} \right) \, ds 
\leq \frac{1}{\tau}  \int\limits_0 ^{kt}  
\left( \int\limits_{\{ u=s \} } d\mathscr{H} _{N-1} (x) \right) \, ds  
\\
 & \leq & 
\frac{f_0 }{\tau g(\tau )} \mathscr{L} ^N (\Omega )  kt, \quad (t\in (0, t_0 ]),
\end{eqnarray*}
and the assertion follows.  
$\hfill \Box $
\\[0.1cm] 
{\sl Proof of Theorem 1: }
The functions
$(u^t -kt)_{+} $ have 
compact support in $\overline{\Omega }$ for all $t\in [0,+\infty ]$, by Lemma 3. 
Hence we derive 
from (\ref{eq1}) the integral identities
\begin{eqnarray}
 0  & = & 
\int\limits_{\Omega } g(|\nabla u|) \frac{\nabla u}{|\nabla u|}  \cdot
\nabla \left( (u^t -kt)_{+} -(u-kt)_{+} \right) \, dx \nonumber \\
 & & -\int\limits_{\Omega } f(|x|,u) 
\left( (u^t -kt)_{+} -(u-kt)_{+} \right) \, 
dx  
\nonumber 
\\
\label{id1}
 & \equiv & I_{1} (t) -I_{2} (t) ,\qquad t\in [0,+\infty ).
\end{eqnarray}
First we claim that 
\begin{equation}
\label{estI2}
I_{2} (t) \geq  o(t).
\end{equation}
To show (\ref{estI2}), we split $I_{2} (t)$ for $0<t\leq u_0 /k$: 
\begin{eqnarray}
I_{2} (t) & = & \int\limits _{ M_1 (t) } \left( f(|x|,u)-f(|x|,kt)\right) 
\left( (u^t -kt)_{+} -(u-kt)_{+} \right) \, dx 
\nonumber 
\\
 & & +\int\limits_{ \Omega } f(|x|,(u-kt)_{+} 
+kt ) \left( (u^t -kt)_{+} -(u-kt)_{+} \right) \, dx 
\nonumber 
\\ 
\label{id2}
 & \equiv & I_{21} (t) +I_{22} (t).
\end{eqnarray} 
By (\ref{ut-u0}) we have
\begin{equation}
\label{utk-uk}
| (u^t (x) -kt )_{+} -(u(x) -kt)_{+} | \leq LRt \qquad \forall x\in \Omega .
\end{equation}
It follows from Lemma 5 and (\ref{utk-uk}), that for $0<t\leq u_0 /k$,
\begin{equation}
\label{estI21}
|I_{21} (t) |  \leq \frac{2}{\tau} f_0 LRc_0  t^2 .
\end{equation}   
Further, in view of (\ref{monot}) we have 
$$
(u^t -kt)_{+} +kt = \left( (u-kt)_{+} +kt \right) ^t .
$$
Hence we obtain by using Remark 3, {\bf 8.},
\begin{eqnarray}
0
 & \leq  & 
\int\limits_{\Omega } \left( F(|x|,(u^t -kt)_{+} +kt) 
-F(|x|,(u-kt)_{+} +kt) \right) \, dx 
\nonumber 
\\
\label{Fest}
 & = & \int\limits_{\Omega } \int\limits_{0} ^{1 }  f(|x|,u^t 
 _{\theta } ) \, d\theta \, 
 \left( (u^t -kt)_{+} -(u-kt)_{+} \right) \, dx ,
\end{eqnarray}
where 
\begin{eqnarray*}
F(r,v) & := & 
\int\limits_{0} ^v f(r,w)\, dw ,\quad (r,v\geq 0),\quad 
\mbox{and } 
\\
u^t _{\theta } & := & 
(1-\theta ) (u-kt)_{+} +\theta (u^t -kt)_{+} 
+kt ,\quad (\theta \in [0,1]).
\end{eqnarray*}
Using 
(\ref{utk-uk}), we obtain, 
\begin{eqnarray*}
 & & |u_{\theta } ^t (x) - (u(x)- kt)_+ -kt | 
  \leq \theta |u^t (x)-u(x)|\leq \theta LRt 
	\\
	 & & \quad \forall x\in \mathbb{R}^N ,\, \theta \in [0,1] .
\end{eqnarray*}
Since the mapping $v\mapsto f(r, v)$ is continuous, uniformly in $r$, this implies
\begin{equation}
\label{f-ft} 
\lim\limits_{t\to 0} \sup \left\{ |f(|x|,u^t _{\theta } (x) ) -f(|x|,(u(x) -kt)_{+} 
 +kt) | :\, 
 x\in \Omega ,\, \theta \in [0,1] \right\} =0.
\end{equation}
Finally,  (\ref{utk-uk}) and (\ref{Fest}) yield, for $0<t\leq u_0 /k$, 
\begin{eqnarray}
I_{22} (t) & \geq  &   \int\limits_{\Omega } \left\{ f(|x|,(u-kt)_{+} +kt) 
-\int\limits_{0}^1  f(|x|,u^t _{\theta }) \, d\theta \right\} 
 \left( (u^t -kt)_{+} -(u-kt)_{+} \right) \, dx   
\nonumber 
\\
 & \geq & 
 -\sup \left\{ \left| f(|x|,u^t _{\theta } ) -f(|x|,(u -kt)_{+} 
 +kt) \right| :\, 
 x\in \Omega ,\, \theta \in [0,1] \right\} \cdot LRt \cdot \mathscr{L} ^N (\Omega ) .
\label{estI22}
\end{eqnarray}
Now  (\ref{estI2}) follows from (\ref{estI21}), (\ref{f-ft}) and (\ref{estI22}).
\\
Next we estimate $I_{1} (t)$.
Let 
\begin{eqnarray*}
G(z) & := & \int_0^z g(s)\, ds 
\quad \mbox{ and }
\\
h(z) & := & G(z) -zg(z), \quad (z\geq 0 ). 
\end{eqnarray*}
Note that $h$ is nonincreasing.
Since $G$ is convex, we have for $t\in (0, u_0 /k ] $,
\begin{eqnarray}
\nonumber
I_1 (t) 
& = & 
\int\limits_{ \Omega \setminus M_2 (t)  } 
g(|\nabla u|) \frac{\nabla u \cdot \nabla u^t }{|\nabla u |} \, dx 
  - \int\limits_{ \Omega \setminus M_1 (t)   } 
	g(|\nabla u|) |\nabla u | \, dx 
\\
	\nonumber
& \leq & 
	\int\limits_{ \Omega \setminus M_2 (t)  } g(|\nabla u|) |\nabla u ^t | \, dx - 
  \int\limits_{ \Omega \setminus M_1 (t)   } g(|\nabla u|) |\nabla u | \, dx
\\
	\nonumber 
& \leq & 
	\int\limits_{ \Omega \setminus M_2 (t)  } \left( G(|\nabla u^t |) - 
	G(|\nabla u|) + g(|\nabla u| ) |\nabla u| \right) \, dx 
\\
	\nonumber 
& & 
- \int\limits_{ \Omega \setminus M_1 (t)  } g(|\nabla u|) |\nabla u | \, dx 
\\
	\nonumber
& = & 
\int\limits_{\Omega \setminus M_2 (t)  } G(|\nabla u^t | )\, dx  - 
\int\limits_{ \Omega \setminus M_1 (t)  } G(|\nabla u | )\, dx 
\\
	\nonumber 
&  & + \int\limits_{ M_2 (t) } h(|\nabla u| )\, dx - 
\int\limits_{ M_1 (t) } h(|\nabla u| )\, dx 
\\
& =: & 
I_{11} (t) - I_{12} (t) + I_{13} (t) - I_{14} (t), \qquad t\in [0, +\infty ).
\label{I1split}
\end{eqnarray}
By Lemma 4 and (\ref{bdry2nhood}) we have
\begin{equation}
\label{nablauest}
\lim\limits_{t\to 0} \sup \left\{ \left| |\nabla u(x) |-\lambda (|x|)\right| :\, 
x\in M_1 (t) \cup M_2 (t) 
 \right\} =0.
\end{equation}
Furthermore, since $\lambda $ is nondecreasing and since $h(z):= 
G(z)-zg(z)$, ($z\in [0, + \infty ) $), is nonincreasing, the function 
$$
p(x):= h\left( \lambda (|x|)\right) , \quad (x\in \mathbb{R}^N ),
$$
satisfies $p=p^{*} =p^t $ $\ \forall 
t\in [0,+\infty ]$. Applying the Hardy-Littlewood inequality 
(\ref{hardylit}) 
we obtain for $t\in (0, u_0 /k]$,
\begin{eqnarray}
 & & \int\limits_{ \Omega \setminus M_1 (t)  } p  \, dx = 
\int\limits_{\mathbb{R} ^N} p\chi (\Omega \setminus M_1 (t)  )  \, dx 
\nonumber 
\\
\label{harlit2} 
 & \leq & 
\int\limits_{\mathbb{R} ^N} p\chi \left( \Omega \setminus M_2 (t)  \right) \, dx = 
\int\limits_{ \Omega \setminus M_2 (t)  } p \, dx .    
\end{eqnarray}
In view of 
(\ref{M1=M2}), 
Lemma 5, (\ref{nablauest}) and (\ref{harlit2}) we obtain
for $t\in (0, u_0 /k ] $,
\begin{eqnarray}
 & & I_{12} (t) -I_{13} (t) = -\int\limits_{ M_1 (t)  } p\,  dx
+\int\limits_{ M_2 (t) } p\,  dx 
\nonumber 
\\
 & & + \int\limits_{ M_1 (t) } \left( p(x) -h\left( |\nabla u|\right) \right) \, dx  
-\int\limits_{ M_2 (t) } \left( p(x) -h\left( |\nabla u|\right) \right) \, dx  
\nonumber 
\\
 & \leq & 
2\mathscr{L}^N \left( M_1 (t)  \right)   
 \sup \left\{ \left| p(x) - 
h\left( |\nabla u(x)| \right) \right| :\, 
x\in M_1 (t) \cup M_2 (t)  \right\}  
\nonumber 
\\
 & \leq & 2c_0 t \sup \left\{ \left| p(x) - 
h\left( |\nabla u(x)| \right) \right| :\, 
x\in M_1 (t) \cup M_2 (t)  \right\}  
 =
 o(t).
\label{bdryest}
\end{eqnarray}
In conclusion, we have by (\ref{id2}),(\ref{estI2}),(\ref{I1split}) and 
(\ref{bdryest})  for $t\in [0, u_0 /k ]$,
\begin{equation}
\label{finalest1}
\int\limits_{\Omega } \left( G(|\nabla (u^t 
-kt)_{+} |) -G(|\nabla (u-kt)_{+} |) \right) \, dx 
\geq o(t). 
\end{equation} 
Now fix $\varepsilon \in(0, u_0 ]$.
Setting $w(x):= \min \{ (u(x)-kt)_{+} ; 
\varepsilon -kt\} $ for $t\in (0,u_0 /k ]$ and $x\in \mathbb{R}^N $,   
we have by (\ref{polyasz1}),
\begin{equation}
\label{finalest2}
\int\limits_{\Omega } G(|\nabla w^t |)\, dx 
\leq  
\int\limits_{\Omega } G(|\nabla w |)\, dx. 
\end{equation}
Since $w^t (x)= \min \{ (u^t (x) -kt)_{+} ; \varepsilon -kt\} $, ($x\in \mathbb{R} ^N $), 
(\ref{finalest1}) together with (\ref{finalest2}) gives for $t\in (0,u_0 /k]$, 
\begin{eqnarray}
 & & \int\limits_{\Omega } \Big( G(|\nabla (u^t 
-\varepsilon )_{+} |) -G(|\nabla (u-\varepsilon )_{+} |) \Big) \, dx  
\nonumber \\
\nonumber
 & = & \int\limits_{\Omega }  \left( G(|\nabla w |) - G(|\nabla w^t |)\right) \, dx
 \\
 \nonumber 
  & &  
+ \int\limits_{\Omega } \Big( G(|\nabla (u^t 
-kt)_{+} |) -G(|\nabla (u-kt)_{+} |) \Big) \, dx 
\\
 & \geq & 
\int\limits_{\Omega } \Big( G(|\nabla (u^t 
-kt)_{+} |) -G(|\nabla (u-kt)_{+} |) \Big) \, dx \geq o(t). 
\label{finalest3}
\end{eqnarray}
In view of Lemma 1, $(u -{\varepsilon } )_{+} $ is locally symmetric in 
direction $x_{1} $. Since the same estimate 
(\ref{finalest3}) can be obtained for CStS in arbitrary directions, Lemma 2 tells us that 
$\{ u>\varepsilon \} $ is an at most countable union of mutually disjoint open balls and 
$(u-\varepsilon )_{+} $ satisfies the symmetry 
property (\ref{decomp1}) of Theorem 1.   Since $\varepsilon $ was arbitrary,
$u$ is locally symmetric in every direction, too. 
Moreover, since $\Omega = \bigcup _{\varepsilon >0 } \{ u>\varepsilon \} $, and $\Omega $ is connected, it must be a ball. The Theorem is proved.
$\hfill \Box $
\\[0.1cm]
{\bf Remark 4.}  It is often possible to derive the 
radial symmetry of the solution of (\ref{eq1}),(\ref{bdry1}) from their  
local symmetry by using other well-known tools 
(see \cite{Br2}, section 10, and \cite{Br3}).
\\
For instance, in the case of the $p$-Laplacian, i.e. if 
$g(z)=z^{p-1} $ for some  $p>1$, $u$ is radially symmetric if 
$f(|x|,\cdot)$ 
satisfies some growth conditions in neighbourhoods of its zero points 
(see \cite{Br3}, Theorem 1).
\\
We mention three typical situations for a general $g$:
\\[0.1cm]
{\bf Theorem 2:} {\sl 
Let $\Omega ,f,g,\lambda $ and $u$ be as in Theorem 1, and  
suppose that one of the following conditions {\bf (a)},{\bf 
(c)} is satisfied:
\\ 
{\bf (a)} $f$ is nonnegative,
\\ 
{\bf (b)} the mapping $r\mapsto f(r,v) $ is strictly decreasing,
\\
{\bf (c)} $f$ is independent of $x$ and the mapping $w\mapsto f(w)$ is nonincreasing.
\\
Then $u$ is radially symmetric and radially decreasing, i.e.  
$\Omega =B_{R} (y)$ for some $ R>0$ and $y\in \mathbb{R}^N $, and 
\begin{equation}
\label{radial}
u=u(r),\quad  \frac{\partial u}{\partial r} \leq 0 ,\qquad ( r=|x-y| ). 
\end{equation}
Moreover, we have $y=0$ in case {\bf (b)}.
}
\\[0.1cm]
{\sl Proof:\/} We use the notations of (\ref{decomp1}).
\\
{\bf (a)} If $r_k >0 $ for some $k\in \{ 1, \ldots ,m\} $, we have by Green's Theorem,
$$
\int\limits_{B_{r_{k}} (y^{k})} f(|x|,u)\, dx =\int\limits_{\partial 
B_{r_{k}} (y^k ) } g(|\nabla u|) \, d\mathscr{H} _{N-1 } (x)  =0, 
$$
which means that $f\equiv 0$ and hence $u=\mbox{const}$ 
in $B_{r_{k}} (y^{k} )$.  
Since $u$ is positive, this implies that we must have $m=1$. Moreover, if $r_1 >0$, then we must have $u=\mbox{const}$ in $B_{r_{1}} (y^{1})$.
\\
{\bf (b)} If the mapping $r\mapsto f(r,v)$ is strictly decreasing, then  
(\ref{eq1}) shows that we must have  
$y^{k} =0$, ($k=1,\ldots ,m$), which proves (\ref{radial}), with $y=0$.
\\
{\bf (c)} Let $\Sigma $ be an arbitrary  $(N-1)$-hyperplane containing the centre $y$ of the ball $\Omega $. Denote by $H$ one of  the two open halfspaces into which $\mathbb{R} ^N $ is split by $\Sigma $, 
and let  $\sigma $ denote reflection in $\Sigma $. Define 
$$
v(x) := u(\sigma x) , \quad (x\in \overline{\Omega \cap H } ).
$$
Since 
$$
\left( g(|y| ) \frac{y}{|y|} - g(|z|) \frac{z}{|z|} \right) \cdot (y-z) \geq 0
$$
for all vectors $y,z\in \mathbb{R} ^N $,
and since the mapping 
$w\mapsto f(w)$ is
nonincreasing, we obtain
\begin{eqnarray*}
0 & \leq & \int\limits_{\Omega \cap H } \left( g(|\nabla u|) \frac{\nabla u}{|\nabla u|} 
-g(|\nabla v|) \frac{\nabla v}{|\nabla v|} \right) \cdot \nabla (u-v) \, dx 
\\
 & = & \int\limits_{\Omega \cap H } ( f( u) - f(  v) )(u-v) \, dx 
 \leq  0.
\end{eqnarray*}
Hence $u= v $ in $\Omega \cap H $. Since $\Sigma $ was arbitrary, the assertion follows. 
$\hfill \Box $ 
\\[0.1cm]
{\bf Acknowledgement:}  I thank  L. 
Damascelli, A. Henrot, W. 
Reichel and 
A. Wagner for helpful discussions.
\\[0.1cm]

\end{document}